\documentclass[12pt]{article}

\usepackage[utf8]{inputenc}
\usepackage[english]{babel}
\usepackage{amsmath,amssymb,amsfonts,amsthm}
\usepackage{bm}
\usepackage{multirow}
\usepackage{xcolor}
\usepackage{hyperref}
\hypersetup{colorlinks=true,linkcolor=blue,filecolor=blue,citecolor = blue,urlcolor=blue} 

\usepackage{graphicx}
\usepackage[scale=0.80]{geometry}

\usepackage{cases}

\theoremstyle{plain}
\newtheorem{theorem}{Theorem}

\theoremstyle{definition}

\theoremstyle{remark}

\newcommand{\NN}{\mathbb{N}}
\newcommand{\ZZ}{\mathbb{Z}}

\newcommand{\prob}[1]{\mathbb{P}\left( #1 \right)}
\newcommand{\esp}[1]{\mathbb{E}\left[ #1 \right]}

\newcommand{\known}{*}
\newcommand{\ci}{?}

\title{Ergodicity of some probabilistic cellular automata\\
  with two letters alphabet via random walks}
\author{Jérôme Casse\\{\small Université Paris-Saclay, CNRS}\\{\small Laboratoire de mathématiques d’Orsay}\\{\small 91405 Orsay, France}}
\date{}

\begin{document}
\maketitle

\begin{abstract}
  Ergodicity of probabilistic cellular automata is a very important issue in the PCA theory. In particular, the question about the ergodicity of all PCA with two-size neighbourhood, two letters alphabet and positive rates is still open. In this article, we do not try to improve this issue, but we show a new kind of proof (to the best knowledge of the author) about the ergodicity of some of those PCA, including also some CA with errors. The proof is based on the study of the boundaries of islands where the PCA is totally decorrelated from its initial condition. The behaviours of these boundaries are the ones of random walks.
\end{abstract}



\section{Intro}

\paragraph{Probabilistic cellular automata.}
In this article, we focus on probabilistic cellular automata with two-size neighbourhood and two-size alphabet $\{0,1\}$.\par

\begin{figure}
  \begin{center}
    \includegraphics{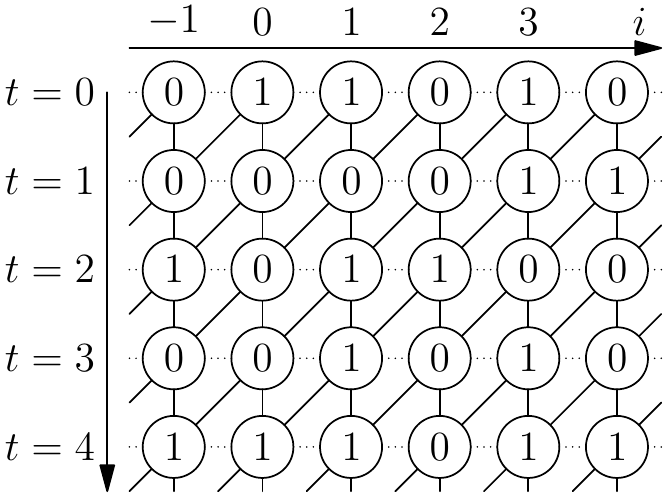}
  \end{center}
  \caption{A realisation of the PCA with parameter $(0.8,0.3,0.5,0.6)$.} \label{fig:real}
\end{figure}

Let $(p_{00},p_{01},p_{10},p_{11}) \in [0,1]^4$ be a quadruplet of real numbers between $0$ and $1$. From this quadruplet and an initial condition $x_0=(x_{0,i})_{i \in \ZZ} \in \{0,1\}^\ZZ$, we can define a Markov chain $X(t)=(X_{i}(t))_{i \in \ZZ}$ on $\{0,1\}^\ZZ$ in such a way :
\begin{itemize}
\item $X(0) = x_0$ and
\item for any $i \in \ZZ$ and any $t \in \NN$,
  \begin{equation}
    X_{i}(t+1) = \begin{cases}
      0 & \text{with probability } 1-p_{X_{i}(t)X_{i+1}(t)},\\
      1 & \text{with probability } p_{X_{i}(t) X_{i+1}(t)}.
    \end{cases}
  \end{equation}
  Moreover, the random variables $(X_{i}(t+1))_{i \in \ZZ}$ are independent knowing $X(t)$.
\end{itemize}
These Markov chains are called \emph{Probabilistic Cellular Automata} (PCA) with two-size neighbourhood and two-size alphabet. In the rest of this article because no confusion is possible, we just called them PCA and the quadruplet $(p_{00},p_{01},p_{10},p_{11})$ the parameter of the PCA. In Figure~\ref{fig:real}, we draw a realisation of one of them.

Even if it is a simple model of Markov chain on uncountable set, it is a very rich model. In particular, let us mention that finding an explicit formula for their invariant measure is a very complex problem solved only when the invariant measure is Markovian~\cite{BGM69,MM14a,CM15,Casse16}.\par
Moreover, knowing if all positive rate PCA with two-size neighbourhood and two-size alphabet\footnote{A PCA with parameter $(p_{00},p_{01},p_{10},p_{11}) \in [0,1]^4$ is called positive rates if $(p_{00}, p_{01},p_{10},p_{11}) \in (0,1)^4$.} are ergodic is still an open problem. It has been solved for many of them, but it still open for around $10\%$ of them\footnote{The $10\%$ corresponds to the volume of the set $\{p \in (0,1)^4 : \text{PCA of parameter } p \text{ is not known to be ergodic}\}$ according to the volume of the hypercube $(0,1)^4$.}. To obtain, the fact that $90\%$ of positive rate PCA with two-size neighbourhood and two-size alphabet are ergodic many techniques more or less complex have been used: coupling~\cite{TOOM90,BMM13,MST19}, cluster expansions~\cite{TOOM90}, contracting maps~\cite{Vasilyev78,TOOM90,MM14a,CM20}, entropy~\cite{KV80,Yaguchi00,DPL02,MST19}, Fourier analysis~\cite{TOOM90,MST19}, weight functions~\cite{HMM19,BKPR22}. In this article, we introduce a new technique based on random walks. If the size of the neighbourhood and the size of the alphabet are sufficiently large, Gàcs proved that there exists a non ergodic positive rates PCA in a two hundreds pages article~\cite{Gacs01}.\par
In addition, they have many applications in combinatorics and statistical physics, see for instance~\cite{Dhar82,BM07,BM98,MM14b,Casse18,CM20,HMM19,SST22} and references therein.

\paragraph{Deterministic cellular automata (with errors).}
When the parameter of a PCA is in $\{0,1\}^4$, PCA are no more probabilistic, but deterministic. They are $16$ of them called \emph{Cellular Automata} (CA). They also have been studied a lot in a various context, see for instance~\cite{Wolfram86,TM87,CD98,Kari05,CC10,Garzon12} and references therein. In the following, the CA with parameter $(p_{00},p_{01},p_{10},p_{11})$ is denoted CA $p_{00}p_{01}p_{10}p_{11}$ (concatenated word).

We call CA $p_{00}p_{01}p_{10}p_{11}$ with error $\epsilon \in [0,1/2]$, the PCA with parameter $(p'_{00},p'_{01},p'_{10},p'_{11})$ where
\begin{equation}
  p'_{xy} = \begin{cases}
    \epsilon & \text{if } p_{xy}=0,\\
    1-\epsilon & \text{if } p_{xy}=1.
  \end{cases}
\end{equation}
In words, an error occurs in the update of each cell with probability $\epsilon$. All of those CA with errors have been proved to be ergodic for any $\epsilon \in (0,1/2]$~\cite{TOOM90,HMM19,MST19}. The last two proved to be ergodic have been the CA $1000$, and the symmetric one $1110$, with errors. They were proved to be ergodic in~\cite{HMM19} via a technical weight function. We give in Section~\ref{sec:proof1000} a new proof for them based on random walks.

\paragraph{Ergodicity of CA.}
A PCA is called ergodic if there exists a probability measure $\mu_\infty$ on $\{0,1\}^\ZZ$ such that, for any initial measure $\mu_0$ such that $X(0) \sim \mu_0$, the measure $\mu_t$ of $X(t)$ weakly converges to $\mu_\infty$.\par

As already mentioned, the ergodicity of all CA (with two-size neighbourhood and two-size alphabet) with error $\epsilon$ have been proved for any $\epsilon \in (0,1/2]$. In this article, we present an idea that permits to prove it for $14$ about $16$ of them. The two CA with errors, for which the idea does not work, when $\epsilon$ is closed to $0$, are CA $0110$ and CA $1001$.

Hence, in this article, we give an alternative proof of the following theorem.
\begin{theorem} \label{thm:CA}
  For any $p=p_{00}p_{01}p_{10}p_{11} \in \{0,1\}^4 \backslash \{1001,0110\}$ and for any $\epsilon \in (0,1/2]$, the CA~$p$ with error $\epsilon$ is ergodic. 
\end{theorem}
\medskip
In fact, for $12$ of them, the idea is quite simple and can be easily generalised to find a sufficient condition about the ergodicity of PCA.\par

\paragraph{Ergodicity of PCA.}
Before expressing the major theorem for PCA, we introduce some helpful notations. First, for any $x \in \{0,1\}$,
\begin{align*}
  & p^{(0)}_{x} = \min(p_{x0},p_{x1}),\ q^{(0)}_{x} = 1 - \max(p_{x0},p_{x1}) \text{ and } r^{(0)}_{x} = 1 - p^{(0)}_{x}-q^{(0)}_{x};\\
  & p^{(1)}_{x} = \min(p_{0x},p_{1x}),\ q^{(1)}_{x} = 1 - \max(p_{0x},p_{1x}) \text{ and } r^{(1)}_{x} = 1 - p^{(1)}_{x}-q^{(1)}_{x};\\
  & p = \min(p_{00},p_{01},p_{10},p_{11}),\ q = 1 - \max(p_{00},p_{01},p_{10},p_{11}) \text{ and } r = 1 - p - q.
\end{align*}
In words, $p^{(0)}_{x}$ (resp. $q^{(0)}_{x}$) is the minimum of probabilities to get $1$ (resp.\ $0$) knowing that the left parent is $x$ and without any knowledge on the right parent. Now, we define the following notations: for any $i \in \{0,1\}$, for any $x \in \{0,1\}$,
\begin{displaymath}
  P^{(i)}_{x} = r p^{(i)}_{x} + (1-r^{(i)}_{x}) p + r^{(i)}_{x} p^{(1-i)}_{x},\  Q^{(i)}_{x} = r q^{(i)}_{x} + (1-r^{(i)}_{x}) q + r^{(i)}_{x} q^{(1-i)}_{x} \text{ and } R_x = r^{(0)}_{x}r^{(1)}_{x}.
\end{displaymath}
Finally, we define the two quantities $\gamma^{(0)}$ and $\gamma^{(1)}$ in the Table~\ref{table:gamma}.
\begin{table}
  $\bullet$ If $r^{(i)}_{0} \leq r^{(i)}_{1}$,\par
  \begin{center}
    \begin{tabular}{|c|c|c|}
      \hline
      If & $Q^{(i)}_{1} \leq Q^{(i)}_{0}$ & $Q^{(i)}_{1} \geq Q^{(i)}_{0}$  \\
      \hline
      $P^{(i)}_{0} \leq P^{(i)}_{1}$  & $\displaystyle \gamma^{(i)} = \frac{Q^{(i)}_{1}}{1 - \left( Q^{(i)}_0 - Q^{(i)}_1 \right)}$ & $\displaystyle \gamma^{(i)} = \frac{Q^{(i)}_1 P^{(i)}_0 + Q^{(i)}_0 \left( 1-P^{(i)}_1 \right)}{1 - \left( P^{(i)}_1 - P^{(i)}_0 \right)} $  \\
      \hline
      $P^{(i)}_{0} \geq P^{(i)}_{1}$  & $\displaystyle \gamma^{(i)} = \frac{Q^{(i)}_{1}}{1 - \left( Q^{(i)}_0 - Q^{(i)}_1 \right)}$ & $\displaystyle \gamma^{(i)} = \frac{Q^{(i)}_{0} + P^{(i)}_{1} \left( Q^{(i)}_{1}-Q^{(i)}_{0} \right)}{1 - \left( Q^{(i)}_{1}-Q^{(i)}_{0} \right) \left( P^{(i)}_{0}-P^{(i)}_{1} \right)}$ \\
      \hline
    \end{tabular}
  \end{center}
  \bigskip\par
  $\bullet$ If $r^{(i)}_{0} \geq r^{(i)}_{1}$,\par
  \begin{center}
    \begin{tabular}{|c|c|c|}
      \hline
      If & $Q^{(i)}_{1} \leq Q^{(i)}_{0}$ & $Q^{(i)}_{1} \geq Q^{(i)}_{0}$  \\
      \hline
      $P^{(i)}_{0} \leq P^{(i)}_{1}$ & $\displaystyle \gamma^{(i)} = \frac{P^{(i)}_{0}}{1 - \left( P^{(i)}_1 - P^{(i)}_0 \right)}$ & $\displaystyle \gamma^{(i)} = \frac{P^{(i)}_{0}}{1 - \left( P^{(i)}_1 - P^{(i)}_0 \right)}$ \\
      \hline
      $P^{(i)}_{0} \geq P^{(i)}_{1}$ & $\displaystyle \gamma^{(i)} = \frac{P^{(i)}_0 Q^{(i)}_1 + P^{(i)}_1 \left( 1-Q^{(i)}_0 \right)}{1 - \left( Q^{(i)}_0 - Q^{(i)}_1 \right)}$ & $\displaystyle \gamma^{(i)} = \frac{P^{(i)}_{1} + Q^{(i)}_{0} \left( P^{(i)}_{0}-P^{(i)}_{1} \right) }{1 - \left( Q^{(i)}_{1}-Q^{(i)}_{0} \right) \left( P^{(i)}_{0}-P^{(i)}_{1} \right)}$\\
      \hline
    \end{tabular}
  \end{center}
  \caption{The value of $\gamma^{(i)}$ for any $i \in \{0,1\}$ according to several conditions.} \label{table:gamma}
\end{table}

\begin{theorem} \label{thm:PCA}
  For any $(p_{00},p_{01},p_{10},p_{11}) \in (0,1)^4$ such that
  \begin{equation} \label{eq:cond}
    2-r > \min \left( r^{(0)}_0,r^{(0)}_1 \right) + \left(1-\gamma^{(0)}\right) \left| r^{(0)}_0 - r^{(0)}_1 \right| + \min \left( r^{(1)}_0,r^{(1)}_1 \right) + \left(1 - \gamma^{(1)} \right) \left| r^{(1)}_0-r^{(1)}_1 \right|,
  \end{equation}
  the PCA with parameter $(p_{00},p_{01},p_{10},p_{11})$ is ergodic. 
\end{theorem}

This Theorem covers and gives an alternative proof for some PCA already known to be ergodic via other techniques, but not for all of them. Due to the fact that the condition~\eqref{eq:cond} is not linear in the parameter $(p_{00},p_{01},p_{10},p_{11})$, we do not have try to compute the volume of the manifold corresponding to it yet. In particular, the Theorem does not cover the four CA $0110$, $1001$, $1000$ and $1110$ with error $\epsilon$ when $\epsilon$ is closed to $0$. Nevertheless, for the two CA~$1000$ and $1110$ with errors $\epsilon$, the idea is improved in Section~\ref{sec:proof1000} to prove their ergodicity for any $\epsilon \in (0,1/2]$. In contrast, we have no hope of adapting the idea used in this paper to get an alternative proof of the ergodicity of the CA $0110$ and $1001$ with error $\epsilon$ when $\epsilon$ is closed to $0$.

\begin{table}
  \begin{center}
    \begin{tabular}{|c|c|c|c|}
      \hline
      left $\backslash$ right & $0$ & $1$ & $\ci$\\ \hline
      $0$ & {$\begin{cases} 0 & \text{w.p. } 1-p_{00} \\ 1 & \text{w.p. } p_{00} \end{cases}$} & {$\begin{cases} 0 & \text{w.p. }  1-p_{01} \\ 1 & \text{w.p. } p_{01} \end{cases}$} & {$\begin{cases} 0 & \text{w.p. } q^{(0)}_{0} \\ 1 & \text{w.p. } p^{(0)}_{0} \\ \ci & \text{w.p. } r^{(0)}_{0} \end{cases}$} \\ \hline
      $1$ & {$\begin{cases} 0 & \text{w.p. } 1-p_{10} \\ 1 & \text{w.p. } p_{10} \end{cases}$} & {$\begin{cases} 0 & \text{w.p. }  1-p_{11} \\ 1 & \text{w.p. } p_{11} \end{cases}$} & {$\begin{cases} 0 & \text{w.p. } q^{(0)}_{1} \\ 1 & \text{w.p. } p^{(0)}_{1} \\ \ci & \text{w.p. } r^{(0)}_{1} \end{cases}$} \\ \hline
      $\ci$ & {$\begin{cases} 0 & \text{w.p. } q^{(1)}_0 \\ 1 & \text{w.p. } p^{(1)}_{0} \\ \ci & \text{w.p. } r^{(1)}_0 \end{cases}$} & {$\begin{cases} 0 & \text{w.p. } q^{(1)}_1 \\ 1 & \text{w.p. } p^{(1)}_1 \\ \ci & \text{w.p. } r^{(1)}_1 \end{cases}$} & {$\begin{cases} 0 & \text{w.p. } q \\ 1 & \text{w.p. } p \\ \ci & \text{w.p. } r \end{cases}$} \\ \hline
    \end{tabular}
  \end{center}
  \caption{Transitions of the envelope PCA with parameter $(p_{00},p_{01},p_{10},p_{11})$.} \label{tab:trans}
\end{table}

\paragraph{Envelope PCA.}
The notion of \emph{envelope} PCA was introduced first in~\cite{BMM13} to prove ergodicity of some PCA using perfect sampling by coupling in the past. It also has been used in~\cite{MST19} to prove ergodicity of some CA with errors. Envelope PCA are closed to \emph{minoring PCA} introduced by Toom~\cite{TOOM90}.\par

The \emph{envelope PCA} of a PCA with parameter $(p_{00},p_{01},p_{10},p_{11})$ is the PCA with two-size neighbours and with a three-size alphabet $\{0,1,\ci\}$. Its initial state is $X_0 =\ \ci^\ZZ$, i.e.\ only the state~$\ci$ is allowed at time~$0$. The state~$\ci$ must be think in the original PCA as a state that depends on the initial condition. The transitions of the envelope PCA according to the parameter $(p_{00},p_{01},p_{10},p_{11})$ are given in Table~\ref{tab:trans}.

\begin{figure}
  \begin{center}
    \includegraphics{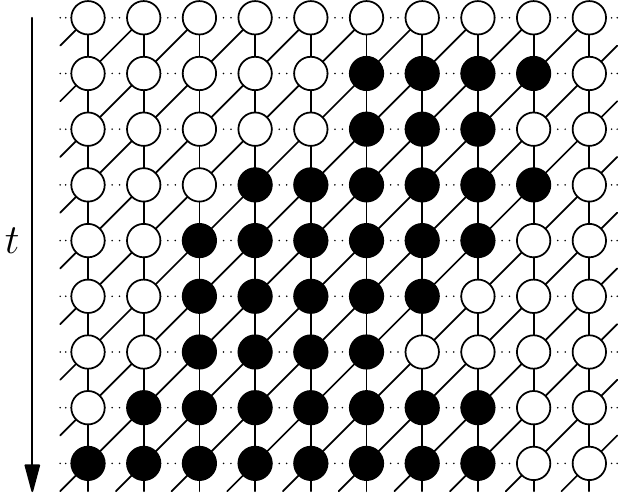} \hspace{2cm} \includegraphics{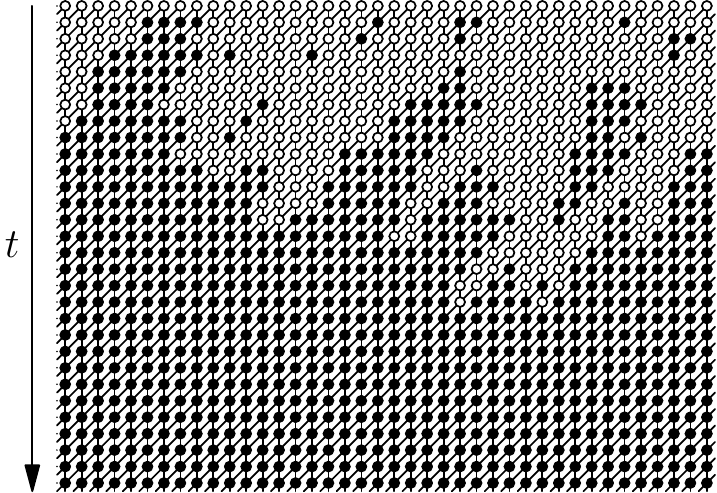}
  \end{center}
  
  \caption{The white cells are cells in state $\ci$ and the black cells are in states $0$ or $1$. On the left, the evolution of a decorrelated island. On the right, the global evolution with multiple decorrelated islands.} \label{fig:ergo}
\end{figure}

Envelope PCA are useful to prove ergodicity of PCA because if states $\ci$ disappear in the envelope PCA then the PCA is ergodic. Now to prove that states $\ci$ disappear, we use the following ideas:
\begin{itemize}
\item With probability $(p+q)^n$, $n$ consecutive cells could go from states $\ci$ to states $0$ or $1$. Such set of consecutive cells is called a \emph{decorrelated island} in the following.
\item In particular, for such decorrelated island, we are interested in the evolution of the locations of its left and right boundaries, denoted $i_{t}$ and $j_{t}$ in $\ZZ$, as well as their states $x_{t}$ and $y_{t}$ in $\{0,1\}$. In particular, the sequence $(i_{t},j_{t},x_{t},y_{t})_{t \geq t_0}$ is a Markov chain where $t_0$ is the creation time of the island.
\item Finally, if $j_{t}-i_{t} \underset{t \to \infty}{\to} \infty$ with a positive probability, then an infinite number of decorrelated islands will grow a.s.\ and the PCA is ergodic.
\end{itemize}
This idea is illustrated on Figure~\ref{fig:ergo}.

\paragraph{Our contribution.}
The main contribution of this paper to the domain is to add the fact that sometimes a geometrical number of cells, decorrelated from the initial conditions, will attach to the boundaries of decorrelated islands. Adding these small contributions permit to go through a recurrent or null recurrent regime for the sizes of the decorrelated islands (i.e.\ $j_{t}-i_{t} \underset{t \to \infty}{\to} 0$ a.s.) to a transient regime for these sizes (i.e.\ $j_{t}-i_{t} \underset{t \to \infty}{\to} \infty$ with positive probability). To the best knowledge of the author, such an idea has not be used to prove ergodicity of PCA in such a way before in the literature.

\paragraph{Content }
In Section~\ref{sec:ergo}, the proof of Theorem~\ref{thm:PCA} is done. In particular, we link the ergodicity of PCA with the transience of random walks describing the sizes of decorrelated islands. The focus is put on the random processes of the left and right boundaries. Then, in Section~\ref{sec:CAerror}, Theorem~\ref{thm:PCA} is applied to twelve CA with errors proving Theorem~\ref{thm:CA} for them. The ergodicity of the two last CA with errors, not proved in Section~\ref{sec:CAerror}, is done in Section~\ref{sec:proof1000}. It is done by some slight improvements helping to describe more precisely the evolution of the boundaries.

\section{Ergodicity of PCA : proof of Theorem~\ref{thm:PCA}} \label{sec:ergo}
Let a PCA of parameter $(p_{00},p_{01},p_{10},p_{11})$. As discuss in the introduction, at time $1$, with probability $(p+q)^{n}$, a decorrelated island of size $n \geq 1$ around the position $0$ is created (in the envelope PCA), i.e.
\begin{displaymath}
  (X_{-\lceil n/2 \rceil + 1}(1),\dots,X_{-1}(1),X_{0}(1),X_{1}(1),\dots,X_{\lfloor n/2 \rfloor}(1)) \in \{0,1\}^n.
\end{displaymath}
This, in fact, can occur at any time $t$. Hence, a.s., there exists $t_0 \geq 1$ such that
\begin{displaymath}
  (X_{-\lceil n/2 \rceil + 1}(t_0),\dots,X_{-1}(t_0),X_{0}(t_0),X_{1}(t_0),\dots,X_{\lfloor n/2 \rfloor}(t_0)) \in \{0,1\}^n.
\end{displaymath}

In the following, we just care about the boundaries of such an island that is:
\begin{itemize}
\item the position of the left border $i_t$ (in the example, it is $i_{t_0} = -\lceil n/2 \rceil + 1$),
\item the position of the right border $j_{t}$ (in the example, it is $j_{t_0} = \lfloor n/2 \rfloor$),
\item the state of the left border $x_{t}$ (in the example, it is $x_{t_0} = X_{-\lceil n/2 \rceil + 1}(t_0)$) and
\item the state of the right border $y_{t}$ (in the example, it is $y_{t_0} = X_{\lfloor n/2 \rfloor}(t_0)$).
\end{itemize}

Indeed, inside $i_{t}$ and $j_{t}-1$, the randomness could be couple to be the same from time $t$ to time $t+1$ and so the states could not be $\ci$'s. Nevertheless, this makes us forgot what is exactly between these two borders and sometimes the exact states of the borders will be forgotten, that is why we need an \emph{additional state $\known$} that means the state is not $\ci$ but we do not remember, with the information kept, if it is a $0$ or a $1$.\par

Let us remark that the evolution of the quadruplet $(i_{t},j_{t},x_{t},y_{t})_{t \geq t_0}$ is a Markov chain and that the random walk $(j_t-i_t)_{t \geq t_0}$ has bounded negative increments (it could decrease of almost $1$). Hence its transcience could be deduced from its asymptotic mean increment.\par

In addition, while $j_{t}-i_{t} \geq 3$, the evolution of the left boundary $(i_{t},x_{t})_{t \geq t_0}$ and the right boundary $(j_{t},y_{t})_{t \geq t_0}$ are independent. In consequence of that, we focus only on the evolution of the right boundary in the next section. The evolution of the left boundary is deduced in the same way.

\subsection{Increment of the right boundary}
Let us suppose that at time $t$, the position of the right boundary is $j_{t} = j$ and its state is $y_{t} = y \in \{0,1,\known\}$.\par

$\bullet$ If $y=0$, then at time $t+1$,
\begin{subnumcases} {\label{eq:t0} (j_{t+1},y_{t+1}) =} 
  (j-1,\known) & w.p. $r^{(1)}_0 r^{(0)}_0$, \label{eq:t0a} \\
  (j-1,0) & w.p. $q^{(1)}_0 r^{(0)}_0$, \label{eq:t0b} \\
  (j-1,1) & w.p. $p^{(1)}_0 r^{(0)}_0$, \label{eq:t0c} \\
  (j,0) & w.p. $q^{(0)}_0 r$, \label{eq:t0d} \\
  (j,1) & w.p. $p^{(0)}_0 r$, \label{eq:t0e} \\
  (j+k+1,0) & w.p. $(1-r^{(0)}_0)(1-r)^{k} q r$, for any $k \in \NN$, \label{eq:t0f} \\
  (j+k+1,1) & w.p. $(1-r^{(0)}_0)(1-r)^{k} p r$, for any $k \in \NN$. \label{eq:t0g}
\end{subnumcases}
These 7 scenarios are illustrated on Figure~\ref{fig:t0}.

\newcommand{\dec}{3.5mm}
\newcommand{\w}{0.24 \textwidth}

\begin{figure}
  \begin{center}
    Initial state\\\smallskip
    \includegraphics[width = 0.5 \textwidth]{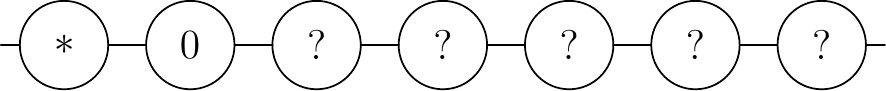}\\\smallskip
    \begin{tabular}{l|l|l}
       \includegraphics[width=\w]{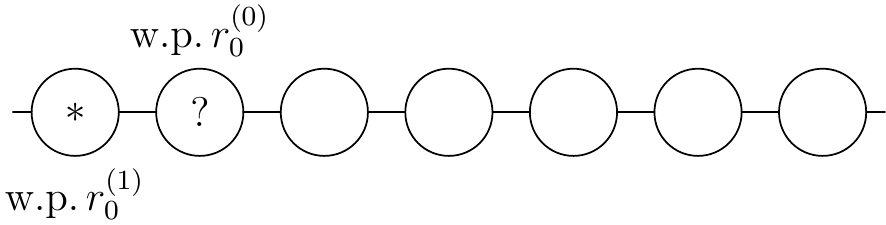}~\raisebox{\dec}{\eqref{eq:t0a}} & & \\
       \includegraphics[width=\w]{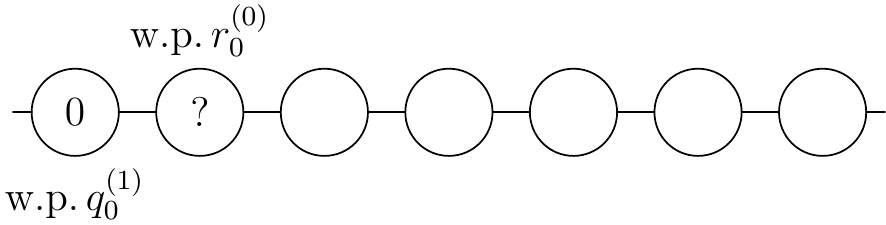}~\raisebox{\dec}{\eqref{eq:t0b}}
     & \includegraphics[width=\w]{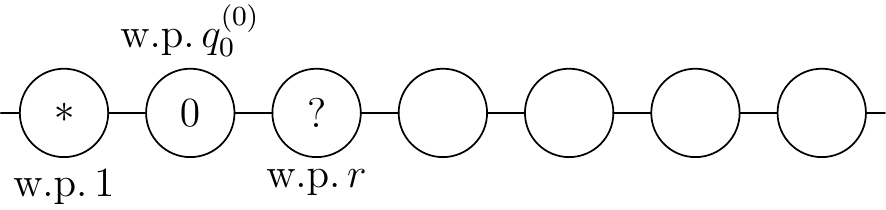}~\raisebox{\dec}{\eqref{eq:t0d}}
     & \includegraphics[width=\w]{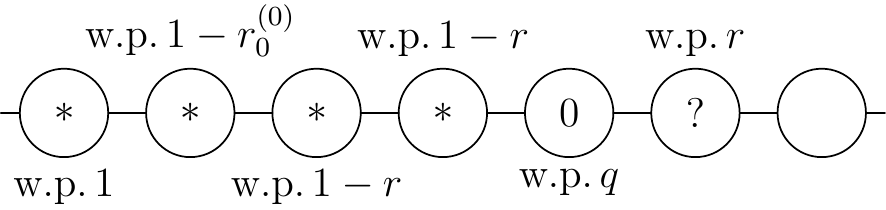}~\raisebox{\dec}{\eqref{eq:t0f}} \\
       \includegraphics[width=\w]{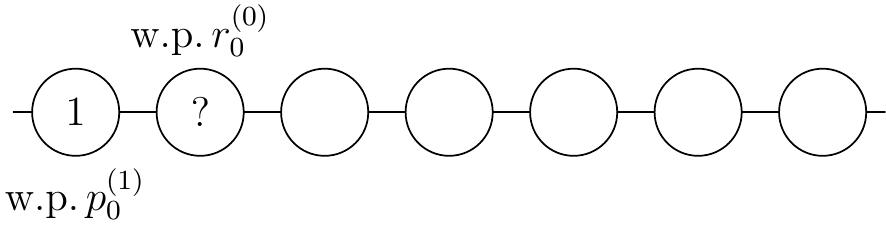}~\raisebox{\dec}{\eqref{eq:t0c}}
     & \includegraphics[width=\w]{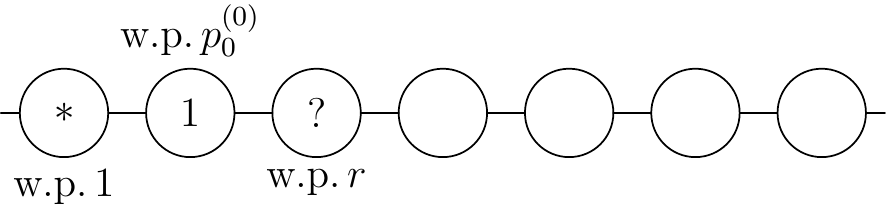}~\raisebox{\dec}{\eqref{eq:t0e}}
     & \includegraphics[width=\w]{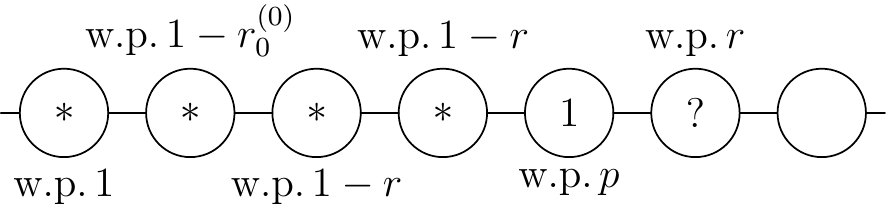}~\raisebox{\dec}{\eqref{eq:t0g}}
    \end{tabular}
  \end{center}
  \caption{In this figure, the initial case is $y = 0$. The 7 possible outcomes are listed below. On each outcome, the probability of each updated cell is written. The third column corresponds to cases~\eqref{eq:t0f} and~\eqref{eq:t0g} with $k=2$.} \label{fig:t0}
\end{figure}

Hence, the law of the right increment $J_{t} = j_{t+1}-j_t$ at time $t$ when the state of the right boundary is $y_{t} = 0$ is
\begin{displaymath}
  \prob{J_{t} = -1~|~y_{t} = 0} = r^{(0)}_0 \text{ and } \prob{J_{t} = k~|~y_{t} = 0} = (1-r^{(0)}_0)(1-r)^kr
\end{displaymath}
whose mean is
\begin{align}
  \esp{J_t~|~y_t=0} & = - r^{(0)}_0 + (1-r^{(0)}_0) r \sum_{k=0}^\infty k (1-r)^k \\
  & = - r^{(0)}_0 + \frac{(1-r^{(0)}_0)(1-r)}{r} = - 1 + \frac{1-r^{(0)}_0}{r}.
\end{align}

$\bullet$ Similarly, if $y=1$, the mean increment is
\begin{displaymath}
  \esp{J_t~|~y_t=1} = - 1 + \frac{1-r^{(0)}_1}{r}.
\end{displaymath}

$\bullet$ Now, if $y = \known$, remember that $\known$ stands for a $0$ or a $1$. Hence, its mean increment is
\begin{align*}
  \esp{J_t~|~y_t=\known} & \geq \min( \esp{J_t~|~y_t=0} ,  \esp{J_t~|~y_t=1} ) \\
                         & = - 1 + \frac{1-\max(r^{(0)}_1,r^{(0)}_0)}{r}.
\end{align*}

\paragraph{The left boundary.}
Similarly and due to the slight asymmetry of the model, the mean increment $I_t$ for the left boundary is
\begin{displaymath}
  \esp{I_t~|~x_t=0} = - \frac{1-r^{(1)}_0}{r},\ \esp{I_t~|~x_t=1} = - \frac{1-r^{(1)}_1}{r} \text{ and }   \esp{I_t~|~x_t=\known} \leq - \frac{1-\max(r^{(1)}_1,r^{(1)}_0)}{r}.
\end{displaymath}

\subsection{Invariant measure of the right boundary state}
\begin{figure}
  \begin{center}
    \includegraphics{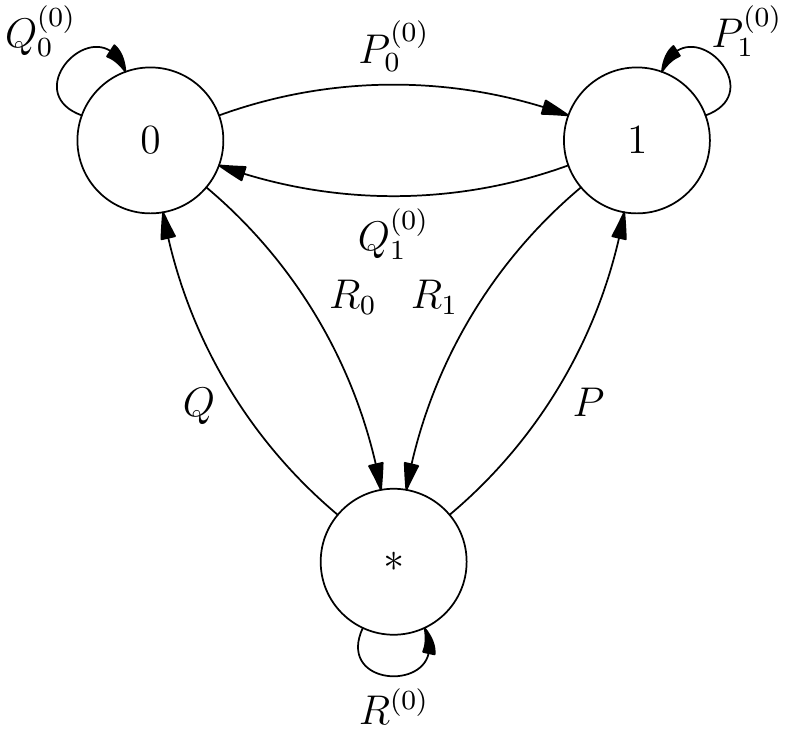}
    \caption{The Markov chain of the state of the boundary.} \label{fig:frontiere}
  \end{center}
\end{figure}

By the previous section, we have seen that there exists a state $w \in \{0,1\}$ whose mean increment is the biggest, the one that satisfies $r^{(0)}_{w} \leq r^{(0)}_{1-w} $. Now, we would like to know how much time the right boundary state spends in the state $w$ compared to the state $1-w$. For that, remark that the right boundary state $(y_t)_{t \geq t_0}$ is a Markov chain on $\{0,1,*\}$ whose transitions are given in Equation~\eqref{eq:t0} and in Figure~\ref{fig:frontiere}. Its transition matrix is
\begin{equation}
  \begin{pmatrix}
    Q^{(0)}_0 & P^{(0)}_0 & R_0 \\
    Q^{(0)}_1 & P^{(0)}_1 & R_1 \\
    Q^{(0)} & P^{(0)} & R^{(0)}
  \end{pmatrix}.
\end{equation}
where
\begin{displaymath}
  Q^{(0)} =  \min \left( Q^{(0)}_0 , Q^{(0)}_1 \right),\ P^{(0)} =  \min \left( P^{(0)}_0 , P^{(0)}_1 \right) \text{, and } R^{(0)} =  1- Q^{(0)} - P^{(0)}.
\end{displaymath}
To determine the transition starting from $*$, as before we remember that $*$ stands for $0$ or $1$, so we just take the coupling between these two that maximises the exact knowledge of $0$ and $1$.

On this Markov chain $(y_t)_{t \geq t_0}$, we are particularly interested in its invariant measure denoted $\nu$ and, specially, by $\nu(w)$. Indeed, by the ergodic theorem and previous section, the asymptotic mean increment of the right boundary is
\begin{equation}
  \lim_{t \to \infty} \esp{J_t} \geq -1 + \frac{1}{r} - \frac{1}{r} \left(\nu(w) r^{(0)}_{w} + (1-\nu(w)) r^{(0)}_{1-w} \right).
\end{equation}
\medskip

To compute $\nu(w)$, three cases can occur:
\begin{itemize}
\item Either ($w=0$ and $Q^{(0)}_1 \leq Q^{(0)}_0$) or ($w=1$ and $P^{(0)}_0 \leq P^{(0)}_1$). In words, that corresponds to the case where the transition from $1-w$ to $w$ and from $*$ to $w$ are equal. In those cases, the invariant measure of $w$ is
  \begin{displaymath}
    \nu(w) =
    \begin{cases}
      \displaystyle \frac{Q^{(0)}_1}{1 - \left( Q^{(0)}_0 - Q^{(0)}_1 \right)} & \text{if } w=0,\\
      \displaystyle \frac{P^{(0)}_0}{1 - \left( P^{(0)}_1 - P^{(0)}_0 \right)} & \text{if } w=1.
    \end{cases}
  \end{displaymath}
  
\item Either ($w=0$ and $Q^{(0)}_1 \geq Q^{(0)}_0$ and $P^{(0)}_0 \leq P^{(0)}_1$) or ($w=1$ and $P^{(0)}_0 \geq P^{(0)}_1$ and $Q^{(0)}_1 \leq Q^{(0)}_0$). In words, that corresponds to the case where the transition from $w$ to $1-w$ and from $*$ to $1-w$ are equal. In those cases, we can deduce, similarly as in the previous case, the invariant measure of $\nu(1-w)$ and so of $\nu(*) + \nu(w)$. The last step is to finish to solve the equations giving the invariant measure to find that
  \begin{displaymath}
    \nu(w) =
    \begin{cases}
      \displaystyle \frac{Q^{(0)}_1 P^{(0)}_0 + Q^{(0)}_0 \left( 1-P^{(0)}_1 \right)}{1-\left( P^{(0)}_1 - P^{(0)}_0 \right)} & \text{if } w=0,\\
      \displaystyle \frac{P^{(0)}_0 Q^{(0)}_1 + P^{(0)}_1 \left( 1-Q^{(0)}_0 \right)}{1 - \left( Q^{(0)}_0 - Q^{(0)}_1 \right)} & \text{if } w=1.
    \end{cases}
  \end{displaymath}
  
\item Either ($Q^{(0)}_1 \geq Q^{(0)}_0$ and $P^{(0)}_0 \geq P^{(0)}_1$). In that case, the system is
  \begin{displaymath}
    \begin{cases}
      \nu(0) = Q^{(0)}_0 \nu(0) + Q^{(0)}_1 \nu(1) + Q^{(0)}_0 \nu(\known) = Q^{(0)}_1 \nu(1) + Q^{(0)}_0 (1-\nu(1)) \\
      \nu(1) = P^{(0)}_0 \nu(0) + P^{(0)}_1 (1-\nu(0)) \\
      \nu(\known) = R_0 \nu(0) + R_1 \nu(1) + R^{(0)} \nu(\known)
    \end{cases}
  \end{displaymath}
  In particular, the two first equations permit to obtain that
  \begin{displaymath}
    \nu(w) =
    \begin{cases}
      \displaystyle \frac{Q^{(0)}_0 + P^{(0)}_1 (Q^{(0)}_1-Q^{(0)}_0)}{1 - (Q^{(0)}_1-Q^{(0)}_0)(P^{(0)}_0-P^{(0)}_1)} & \text{if } w=0,\\
      \displaystyle \frac{P^{(0)}_1 + Q^{(0)}_0 (P^{(0)}_0-P^{(0)}_1) }{1 - (Q^{(0)}_1-Q^{(0)}_0)(P^{(0)}_0-P^{(0)}_1)} & \text{if } w=1.
    \end{cases}
  \end{displaymath}
\end{itemize}

In fact, in all these cases, the value $\nu(w)$ is equal to $\gamma^{(0)}$ as defined in Table~\ref{table:gamma}. Hence, by the ergodic theorem, the asymptotic of the mean increment of the right boundary is
\begin{align*}
  \lim_{t \to \infty} \esp{J_t} \geq & -1 + \frac{1}{r} - \frac{1}{r} \left(\gamma^{(0)} \min(r^{(0)}_0,r^{(0)}_1) + (1-\gamma^{(0)})  \max(r^{(0)}_0,r^{(0)}_1) \right)\\
  = & -1 + \frac{1}{r} - \frac{1}{r} \left(\min (r^{(0)}_0,r^{(0)}_1) + \left( 1-\gamma^{(0)} \right) \left| r^{(0)}_0 - r^{(0)}_1 \right|   \right).
\end{align*}

\paragraph{The left boundary.}
Similarly, the asymptotic of the mean increment of the left boundary is
\begin{displaymath}
  \lim_{t \to \infty} \esp{I_t} \leq - \frac{1}{r} + \frac{1}{r}  \left(\min(r^{(1)}_0,r^{(1)}_1) + \left( 1-\gamma^{(1)} \right) \left| r^{(1)}_0 - r^{(1)}_1 \right| \right).
\end{displaymath}

\subsection{Conclusion}
The sequence $(j_t-i_t)_{t \geq t_0}$ is a random walk with an asymptotic drift
\begin{displaymath}
  D= \lim_{t \to \infty} \esp{J_t} - \lim_{t \to \infty} \esp{I_t}
\end{displaymath}
and bounded negative increments. So, if $D > 0$, it is transient and the PCA is ergodic. But, by the previous section, the drift $D$ is bigger than
\begin{align*}
  -1 + \frac{2}{r} - \frac{1}{r} \left(\min(r^{(0)}_0,r^{(0)}_1) + \left(1-\gamma^{(0)}\right) \left| r^{(0)}_0 - r^{(0)}_1 \right| + \min(r^{(1)}_0,r^{(1)}_1) + \left(1-\gamma^{(1)} \right) \left| r^{(1)}_0 - r^{(1)}_1 \right| \right)
\end{align*}
that is positive if Equation~\eqref{eq:cond} holds.\qed

\section{Applications to CA with errors} \label{sec:CAerror}
In this section, we apply condition~\eqref{eq:cond} to CA with errors $\epsilon$. We distinguish four families for the twelve PCA we focus on. For each family, we show how to apply the condition~\eqref{eq:cond} to one of the CA in the family. For the other CA in the same family, the condition applies similarly.
\begin{itemize}
\item The first family is CA $0000$ and $1111$ for which the condition is trivial due to the fact that $r=0$.
\item The second family is CA $0011$, $0101$, $1010$ and $1100$. 
\item The third family is CA $0001$ and $0111$.
\item The fourth family is CA $0010$, $0100$, $1011$ and $1101$.
\end{itemize}
The condition does not apply for CA $1000$ and CA $1110$ with error $\epsilon$ close to $0$. These two cases are slightly different and are done in Section~\ref{sec:proof1000}.\par
\bigskip
In all this section, let $\epsilon$ be any real number in $(0,1/2]$.

\subsection{CA $\mathbf{0000}$}
This case is trivial. Indeed, the parameters $r=r^{(0)}_{0}=r^{(0)}_{1}=r^{(1)}_{0}=r^{(1)}_{1}=0$, and so the Equation~\eqref{eq:cond} is $2 > 0$.

\subsection{CA $\mathbf{0011}$}
It corresponds to the PCA with parameter $(\epsilon,\epsilon,1-\epsilon,1-\epsilon)$. In particular, the parameters $r^{(0)}_0 = r^{(0)}_1 = 0$ and $r= r^{(1)}_0 = r^{(1)}_1 = 1-2\epsilon$ when $\epsilon \in (0,1/2]$. Equation~\eqref{eq:cond} becomes $1+ 2 \epsilon > 1- 2 \epsilon$ that holds if $\epsilon >0$.

\subsection{CA $\mathbf{0001}$}
It corresponds to the PCA with parameter $(\epsilon,\epsilon,\epsilon,1-\epsilon)$. In that case,
\begin{center}
  \begin{tabular}{lll}
    $p^{(0)}_0 = p^{(1)}_0 = \epsilon$, & $q^{(0)}_0 = q^{(1)}_0 = 1-\epsilon$, &$r^{(0)}_0 = r^{(1)}_0 = 0$;\\
    $p = p^{(0)}_1 = p^{(1)}_1 = \epsilon$, & $q = q^{(0)}_1 = q^{(1)}_1 = \epsilon$, &$r= r^{(0)}_1 = r^{(1)}_1 = 1-2\epsilon$;\\
    \multicolumn{2}{l}{$P^{(0)}_0 = P^{(1)}_0 = P^{(1)}_0 = P^{(1)}_1 = 2 \epsilon (1-\epsilon)$,} \\
    $Q^{(0)}_0 =Q^{(1)}_0 = 1 - 2 \epsilon (1-\epsilon)$, & $Q^{(0)}_1 =Q^{(1)}_1 =2 \epsilon (1-\epsilon)$.
  \end{tabular}
\end{center}
So, for both $i \in \{0,1\}$, $r^{(i)}_0 \leq r^{(i)}_1$ and $Q^{(i)}_1 < Q^{(i)}_0$ when $\epsilon \in (0,1/2]$. Hence, Table~\ref{table:gamma} gives $\gamma^{(0)} = \gamma^{(1)} = 1/2$. Then, Equation~\ref{eq:cond} becomes $1 + 2 \epsilon > 1 - 2 \epsilon$ that holds if $\epsilon > 0$.

\subsection{CA $\mathbf{0010}$}
It corresponds to the PCA with parameter $(\epsilon,\epsilon,1-\epsilon,\epsilon)$. In that case,
\begin{center}
  \begin{tabular}{lll}
    $p^{(0)}_0 = p^{(1)}_1 = \epsilon$, & $q^{(0)}_0 = q^{(1)}_1 = 1-\epsilon$, &$r^{(0)}_0 = r^{(1)}_1 = 0$;\\
    $p = p^{(0)}_1 = p^{(1)}_0 = \epsilon$, & $q = q^{(0)}_1 = q^{(1)}_0 = \epsilon$, &$r= r^{(0)}_1 = r^{(1)}_0 = 1-2\epsilon$;\\
    \multicolumn{2}{l}{$P^{(0)}_0 = P^{(1)}_0 = P^{(1)}_0 = P^{(1)}_1 = 2 \epsilon (1-\epsilon)$,} \\
    \multicolumn{2}{l}{$Q^{(0)}_0 = Q^{(1)}_0 = Q^{(1)}_0 = Q^{(1)}_1 = 1 - 2 \epsilon (1-\epsilon)$.} \\
  \end{tabular}
\end{center}
So, by Table~\ref{table:gamma}, $\gamma^{(0)} = 1-2 \epsilon (1-\epsilon)$ and $\gamma^{(1)} = 2 \epsilon (1-\epsilon)$. Then, Equation~\eqref{eq:cond} becomes $1+2 \epsilon > 1 - 2 \epsilon$ that holds if $\epsilon > 0$.\par

\section{Ergodicity of CA $\mathbf{1000}$ with error $\epsilon \in (0,1/2)$} \label{sec:proof1000}
It corresponds to the PCA with parameter $(1-\epsilon,\epsilon,\epsilon,\epsilon)$.

\subsection{Improvements}
\begin{figure}
  \begin{center}
    \includegraphics{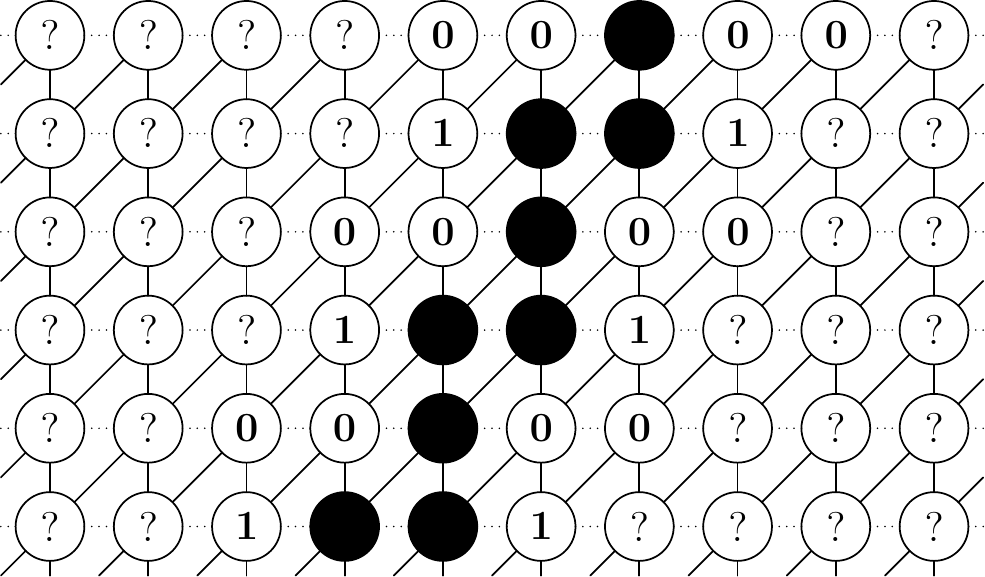}
  \end{center}
  \caption{The evolution of a decorrelated island of CA $1000$ with errors $\epsilon$ when none error occurs.} \label{fig:DC1000}
\end{figure}

The improvement listed below are done to treat the fact that, without error, the state of the boundaries of a decorrelated island oscillate between two consecutive $0$ and one $1$, see~Figure~\ref{fig:DC1000}.

\paragraph{First improvement.}
The first improvement is to consider boundaries of size~$2$ instead of size~$1$. Hence, at a time $t_0$ corresponding to a creation of a decorrelated island (in the envelope PCA),
\begin{displaymath}
  (X_{-\lceil n/2 \rceil + 1}(t_0),\dots,X_{-1}(t_0),X_{0}(t_0),X_{1}(t_0),\dots,X_{\lfloor n/2 \rfloor}(t_0)) \in \{0,1\}^n,
\end{displaymath}
we care, as before, about the positions $i_t$ and $j_t$ of the left and right boundaries (in the example, $i_{t_0}=-\lceil n/2 \rceil + 1$ and $j_{t_0}=\lfloor n/2 \rfloor$) and, now, about 
\begin{itemize}
\item the two states of the left boundary $x_t \in \{0,1,\known\}^2$ : in the example,
  \begin{displaymath}
    x_{t_0} = \left( X_{-\lceil n/2 \rceil + 1}(t_0),X_{t_0,-\lceil n/2 \rceil + 2}(t_0) \right),
  \end{displaymath}
  
\item the two states of the right boundary $y_t \in \{0,1,\known\}^2$ : in the example,
  \begin{displaymath}
    y_{t_0} = \left( X_{\lfloor n/2 \rfloor-1}(t_0),X_{\lfloor n/2 \rfloor}(t_0) \right).
  \end{displaymath}
\end{itemize}
Let us remark that, while $j_t-i_t \geq 5$, the evolution of the quadruplet $(i_t,j_t,x_t,y_t)_{t \geq t_0}$ is a Markov chain and the left and right boundaries are independent.

\paragraph{Second improvement.}
The second improvement is to consider the evolution of two functions $\tilde{i}_t$ and $\tilde{j}_t$ depending both on the positions $i_t$ and $j_t$ and on the states $x_t$ and $y_t$ of the boundaries. These functions are
\begin{align*}
  \tilde{i}_t = \begin{cases}
    i_t & \text{if } x_t \in \{(1,0),(1,1),(1,\known)\},\\
    i_t + 1/2 & \text{if } x_t \in \{(0,0)\},\\
    i_t + 1 & \text{if } x_t \in \{(0,1),(\known,1),(0,\known),(\known,\known),(\known,0)\},\\
  \end{cases}
\end{align*}
and, similarly,
\begin{align*}
  \tilde{j}_t = \begin{cases}
    j_t & \text{if } y_t \in \{(0,1),(1,1),(\known,1)\},\\
    j_t-1/2 & \text{if } y_t \in \{(0,0)\},\\
    j_t-1 & \text{if } y_t \in \{(1,0),(1,\known),(\known,0),(\known,\known),(0,\known)\}.
  \end{cases}
\end{align*}

This improvement is done for a technical reason that simplify the study of the size of the decorrelated island. Indeed, on Figure~\ref{fig:DC1000}, for the drawn decorrelated island, the value of $(\tilde{j}_t-\tilde{i}_t)_t$ is $3$ at any time $t$, whereas the values of $(j_t-i_t)_{t}$ oscillates between $3$ and $4$ according to the parity of $t$. Hence, $\tilde{i}_t$ and $\tilde{j}_t$ reflect more precisely the evolution of the size of the decorrelated island.

In addition, the absolute difference between $\tilde{j}_t-\tilde{i}_t$ and $j_t-i_t$ is bounded by $2$. Hence, if one goes to $\infty$, the other one goes too.

\subsection{Increments of $\mathbf{\tilde{i}_t}$ and $\mathbf{\tilde{j}_t}$ }
\paragraph{Right boundary.}
Let us suppose that at time $t$, $\tilde{j}_t = \tilde{j}$ and $y_t = y$.

\begin{figure}
  \begin{center}
    \includegraphics[width = 0.5 \textwidth]{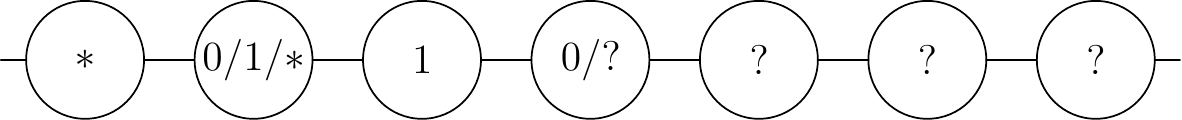}\\\smallskip
    \begin{tabular}{l|l|l}
      \includegraphics[width=\w]{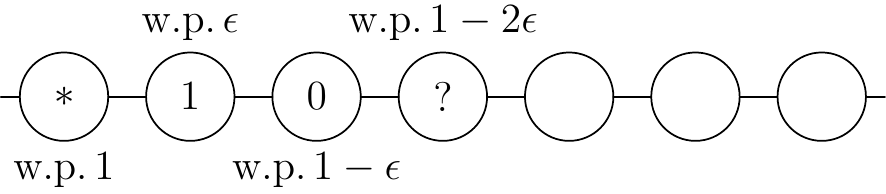}~\raisebox{\dec}{\eqref{eq:t1a}}
    & \includegraphics[width=\w]{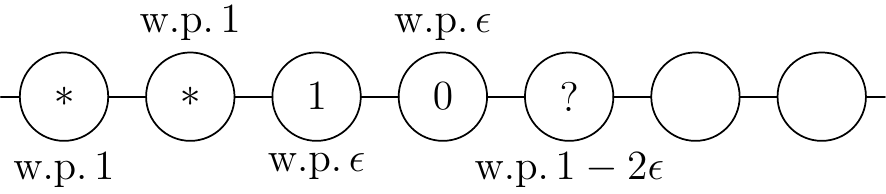}~\raisebox{\dec}{\eqref{eq:t1e}}
    & \includegraphics[width=\w]{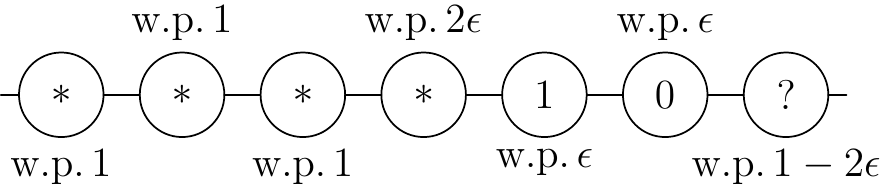}~\raisebox{\dec}{\eqref{eq:t1i}}\\
      \includegraphics[width=\w]{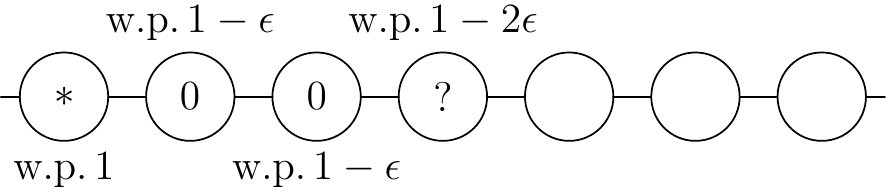}~\raisebox{\dec}{\eqref{eq:t1b}}
    & \includegraphics[width=\w]{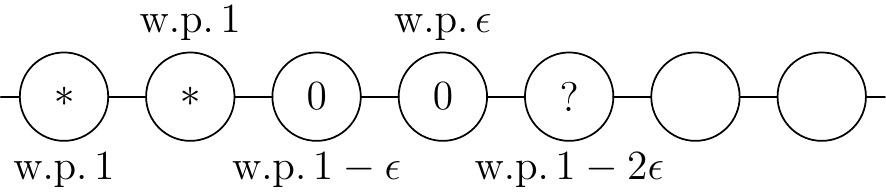}~\raisebox{\dec}{\eqref{eq:t1f}}
    & \includegraphics[width=\w]{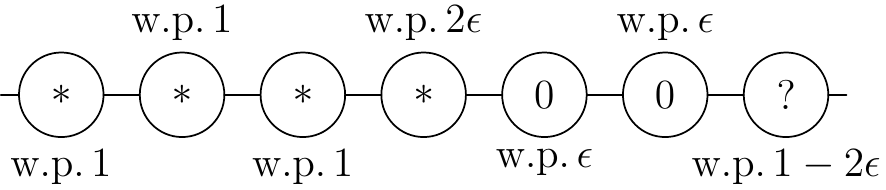}~\raisebox{\dec}{\eqref{eq:t1j}}\\
      \includegraphics[width=\w]{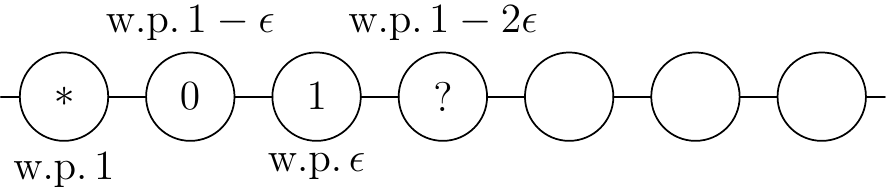}~\raisebox{\dec}{\eqref{eq:t1c}}
    & \includegraphics[width=\w]{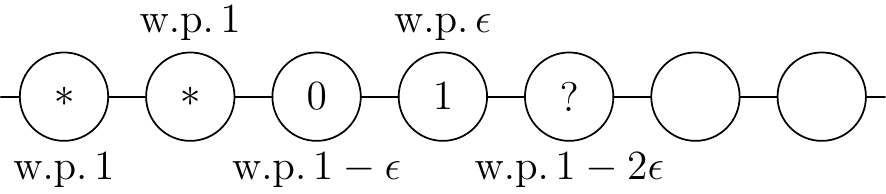}~\raisebox{\dec}{\eqref{eq:t1g}}
    & \includegraphics[width=\w]{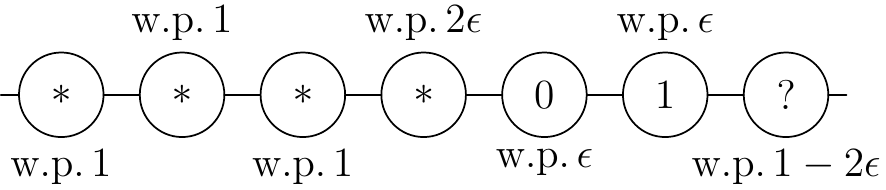}~\raisebox{\dec}{\eqref{eq:t1k}}\\
      \includegraphics[width=\w]{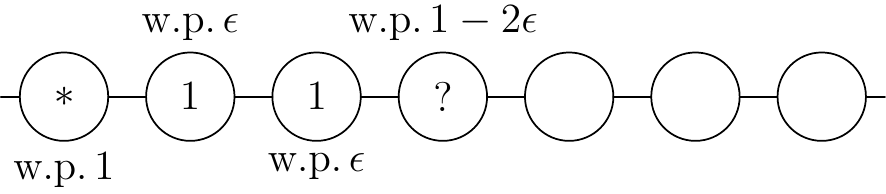}~\raisebox{\dec}{\eqref{eq:t1d}}
    & \includegraphics[width=\w]{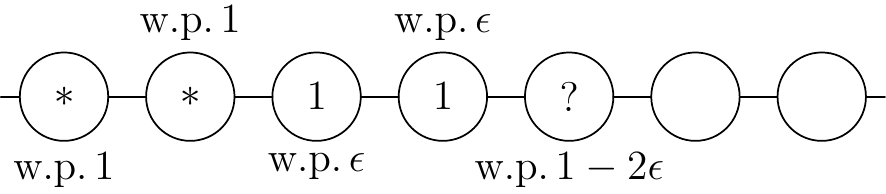}~\raisebox{\dec}{\eqref{eq:t1h}}
    & \includegraphics[width=\w]{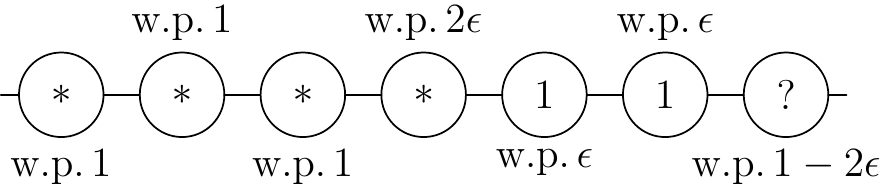}~\raisebox{\dec}{\eqref{eq:t1l}}
    \end{tabular}
  \end{center}
 \caption{In this figure, the initial case on the top corresponds to $y \in \{(0,1),(1,1),(\known,1),(1,0)\}$. The 12 possible outcomes are listed below. On each outcome, the probability of each updated cell is written. The third column corresponds to cases~\eqref{eq:t1i}, \eqref{eq:t1j}, \eqref{eq:t1k} and~\eqref{eq:t1l} with $k=1$.} \label{fig:t1}
\end{figure}

$\bullet$ If $y \in S_1 = \{(0,1),(1,1),(\known,1),(1,0)\}$, then at time $t+1$
\begin{subnumcases} {\label{eq:t1} (\tilde{j}_{t+1},y_{t+1}) =}
  (\tilde{j}-1,(1,0)) & w.p. $\epsilon (1-\epsilon) (1-2\epsilon)$, \label{eq:t1a} \\
  (\tilde{j}-1/2,(0,0)) & w.p. $(1-\epsilon)^2 (1-2\epsilon)$, \label{eq:t1b} \\
  (\tilde{j},(0,1)) & w.p. $(1-\epsilon) \epsilon (1-2\epsilon)$, \label{eq:t1c} \\
  (\tilde{j},(1,1)) & w.p. $\epsilon^2 (1-2\epsilon)$, \label{eq:t1d} \\
  (\tilde{j},(1,0)) & w.p. $\epsilon^2 (1-2\epsilon)$, \label{eq:t1e} \\
  (\tilde{j}+1/2,(0,0)) & w.p. $(1-\epsilon) \epsilon (1-2\epsilon)$, \label{eq:t1f} \\
  (\tilde{j}+1,(0,1)) & w.p. $(1-\epsilon) \epsilon (1-2\epsilon)$, \label{eq:t1g} \\
  (\tilde{j}+1,(1,1)) & w.p. $\epsilon^2 (1-2\epsilon)$, \label{eq:t1h} \\
  (\tilde{j}+k+1,(1,0)) & w.p. $(2 \epsilon)^k \epsilon^2 (1-2\epsilon)$ for any $k \geq 0$, \label{eq:t1i} \\
  (\tilde{j}+k+3/2,(0,0)) & w.p. $(2 \epsilon)^k \epsilon^2 (1-2\epsilon)$ for any $k \geq 0$, \label{eq:t1j} \\
  (\tilde{j}+k+2,(0,1)) & w.p. $(2 \epsilon)^k \epsilon^2 (1-2\epsilon)$ for any $k \geq 0$,\label{eq:t1k} \\
  (\tilde{j}+k+2,(1,1)) & w.p. $(2 \epsilon)^k \epsilon^2 (1-2\epsilon)$ for any $k \geq 0$. \label{eq:t1l} 
\end{subnumcases}
These 12 scenarios are illustrated on Figure~\ref{fig:t1}.

Hence, the mean of the right increment at time $t$ when the state of the right boundary is $y_t \in S_1 = \{(0,1),(1,1),(\known,1),(1,0)\}$ is
\begin{align*}
  \esp{\tilde{J}_t~\big|~y_t \in S_1} & = -\frac{1}{2} + \frac{5\epsilon}{2} + \frac{7 \epsilon^2}{2} + \frac{8\epsilon^3}{1-2\epsilon}
\end{align*}

$\bullet$ If $y = (0,0)$, then at time $t+1$
\begin{subnumcases} {\label{eq:t00} (\tilde{j}_{t+1},y_{t+1}) =}
  (\tilde{j}-3/2,(\known,0)) & w.p. $(1-2\epsilon) \epsilon (1-2\epsilon)$, \label{eq:t00a} \\
  (\tilde{j}-3/2,(1,0)) & w.p. $\epsilon^2 (1-2\epsilon)$, \label{eq:t00b} \\
  (\tilde{j}-1,(0,0)) & w.p. $\epsilon^2 (1-2\epsilon)$, \label{eq:t00c} \\
  (\tilde{j}-1/2,(\known,1)) & w.p. $(1-2\epsilon) (1-\epsilon) (1-2\epsilon)$, \label{eq:t00d} \\
  (\tilde{j}-1/2,(0,1)) & w.p. $\epsilon (1-\epsilon) (1-2\epsilon)$, \label{eq:t00e} \\
  (\tilde{j}-1/2,(1,1)) & w.p. $\epsilon (1-\epsilon) (1-2\epsilon)$, \label{eq:t00f} \\
  (\tilde{j}-1/2,(1,0)) & w.p. $(1-\epsilon) \epsilon (1-2\epsilon)$, \label{eq:t00g} \\
  (\tilde{j},(0,0)) & w.p. $\epsilon^2 (1-2\epsilon)$, \label{eq:t00h} \\
  (\tilde{j}+1/2,(0,1)) & w.p. $\epsilon^2 (1-2\epsilon)$, \label{eq:t00i} \\
  (\tilde{j}+1/2,(1,1)) & w.p. $(1-\epsilon) \epsilon (1-2\epsilon)$, \label{eq:t00j} \\
  (\tilde{j}+k+1/2,(1,0)) & w.p. $(2 \epsilon)^k \epsilon^2 (1-2\epsilon)$ for any $k \geq 0$, \label{eq:t00k} \\
  (\tilde{j}+k+1,(0,0)) & w.p. $(2 \epsilon)^k \epsilon^2 (1-2\epsilon)$ for any $k \geq 0$, \label{eq:t00l} \\
  (\tilde{j}+k+3/2,(0,1)) & w.p. $(2 \epsilon)^k \epsilon^2 (1-2\epsilon)$ for any $k \geq 0$, \label{eq:t00m} \\
  (\tilde{j}+k+3/2,(1,1)) & w.p. $(2 \epsilon)^k \epsilon^2 (1-2\epsilon)$ for any $k \geq 0$. \label{eq:t00n}
\end{subnumcases}
  
\begin{figure}
  \begin{center}
    \includegraphics[width = 0.5 \textwidth]{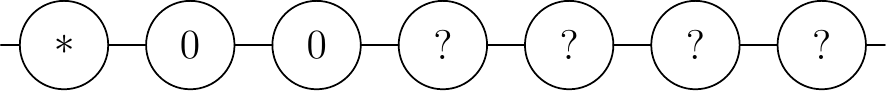}\\\smallskip
    \begin{tabular}{l|l|l}
      \includegraphics[width=\w]{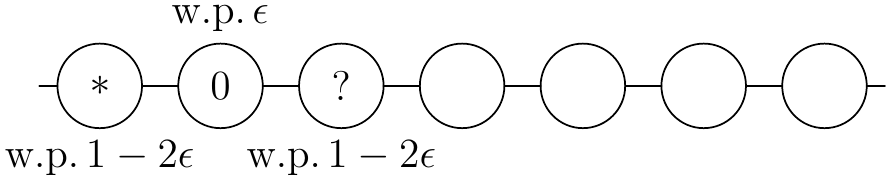}~\raisebox{\dec}{\eqref{eq:t00a}} & & \\
      \includegraphics[width=\w]{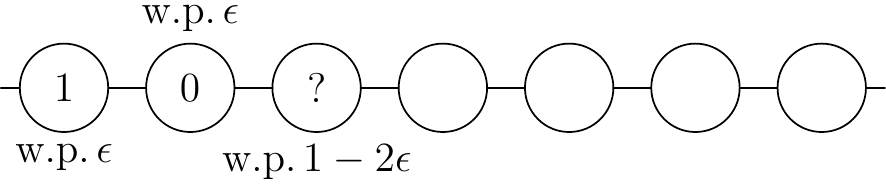}~\raisebox{\dec}{\eqref{eq:t00b}}
    & \includegraphics[width=\w]{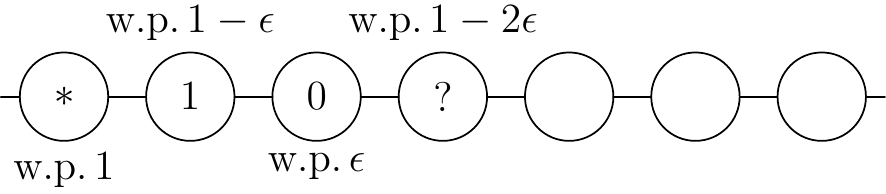}~\raisebox{\dec}{\eqref{eq:t00g}}
    & \includegraphics[width=\w]{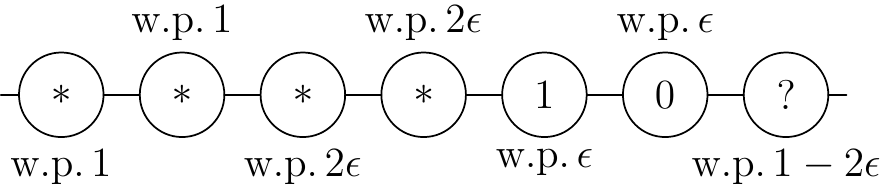}~\raisebox{\dec}{\eqref{eq:t00k}}\\
      \includegraphics[width=\w]{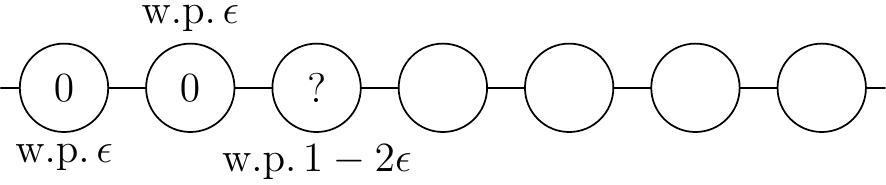}~\raisebox{\dec}{\eqref{eq:t00h}}
    & \includegraphics[width=\w]{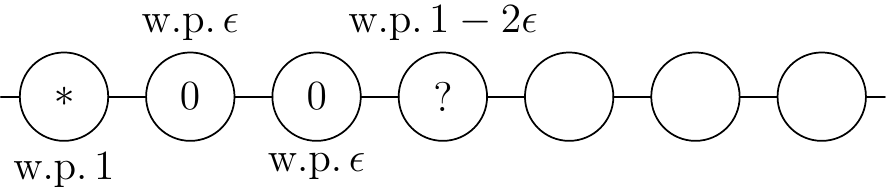}~\raisebox{\dec}{\eqref{eq:t00f}}
    & \includegraphics[width=\w]{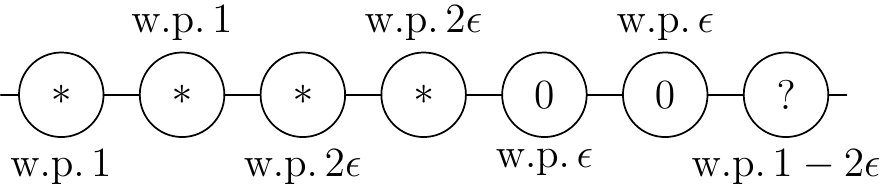}~\raisebox{\dec}{\eqref{eq:t00l}}\\
      \includegraphics[width=\w]{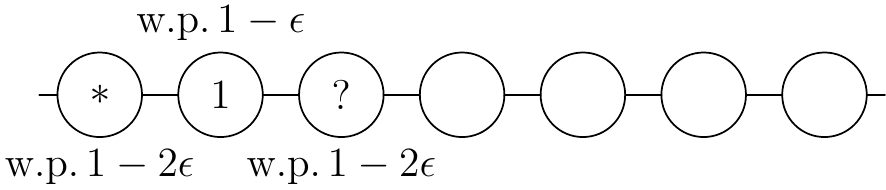}~\raisebox{\dec}{\eqref{eq:t00d}} & & \\
      \includegraphics[width=\w]{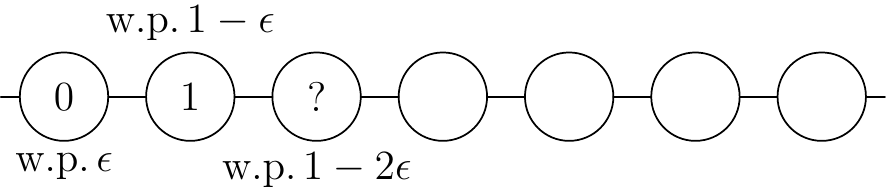}~\raisebox{\dec}{\eqref{eq:t00e}}
    & \includegraphics[width=\w]{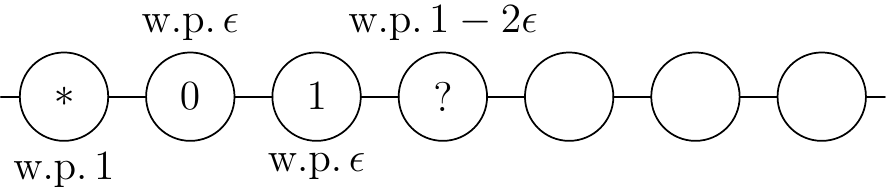}~\raisebox{\dec}{\eqref{eq:t00i}}
    & \includegraphics[width=\w]{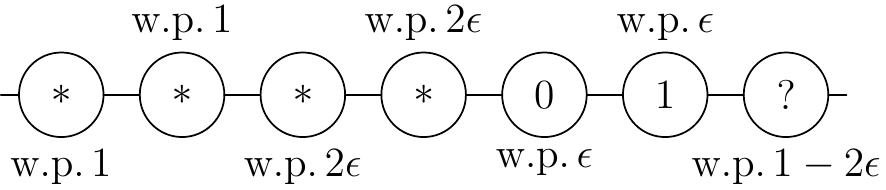}~\raisebox{\dec}{\eqref{eq:t00m}}\\
      \includegraphics[width=\w]{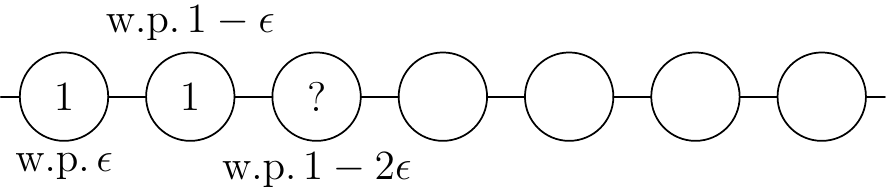}~\raisebox{\dec}{\eqref{eq:t00f}}
    & \includegraphics[width=\w]{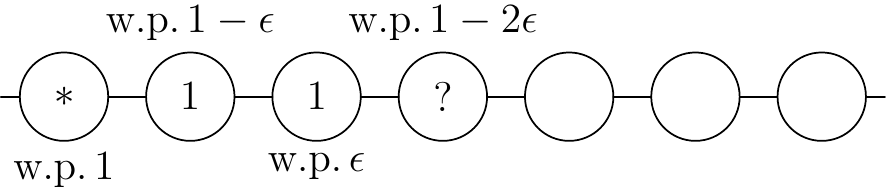}~\raisebox{\dec}{\eqref{eq:t00j}}
    & \includegraphics[width=\w]{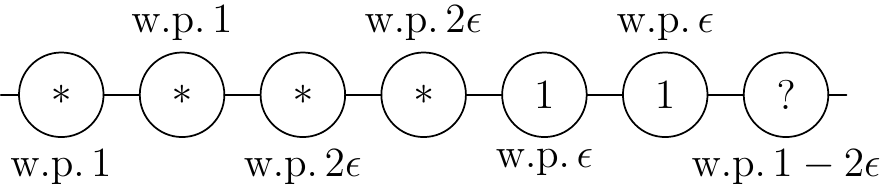}~\raisebox{\dec}{\eqref{eq:t00n}}
    \end{tabular}
  \end{center}
  \caption{In this figure, the initial case on the top corresponds to $y = (0,0)$. The $14$ possible outcomes are listed below. On each outcome, the probability of each updated cell is written. The third column corresponds to cases~\eqref{eq:t00k}, \eqref{eq:t00l}, \eqref{eq:t00m} and~\eqref{eq:t00n} with $k=2$.} \label{fig:t00}
\end{figure}

These 14 scenarios are illustrated on Figure~\ref{fig:t00}.

Hence, the mean increment at time $t$ when the right boundary is $y_t = (0,0)$ is
\begin{displaymath}
  \esp{\tilde{J}_t~\big|~y_t = (0,0)} = -\frac{1}{2} + \frac{15 \epsilon^2}{2} + 6 \epsilon^3 + \frac{16\epsilon^4}{1-2\epsilon}.
\end{displaymath}

$\bullet$ At this step, we should look for the last case: if $y_t = (\known,0)$. In fact, that case corresponds to either $(0,0)$ either $(1,0)$, so its mean increment is, at least, for any $\epsilon \in (0,1/2)$,
\begin{align*}
  \esp{\tilde{J}_t~\big|~y_t = (\known,0)} & \geq \min\left( \esp{\tilde{J}_t~\big|~y_t = (0,0)},\esp{\tilde{J}_t~\big|~y_t = (1,0)} \right)\\
                                           & = \esp{\tilde{J}_t~\big|~y_t = (0,0)}.
\end{align*}

Hence, the asymptotic mean increment of the right boundary is
\begin{displaymath}
  \lim_{t \to \infty} \esp{\tilde{J}_t} \geq -\frac{1}{2} + \frac{15 \epsilon^2}{2} + 6 \epsilon^3 + \frac{16\epsilon^4}{1-2\epsilon}.
\end{displaymath}

\paragraph{Left boundary.}
Similarly, the asymptotic mean increment of the left boundary is
\begin{displaymath}
  \lim_{t \to \infty} \esp{\tilde{I}_t} \leq -\frac{1}{2} - \frac{15 \epsilon^2}{2} - 6 \epsilon^3 - \frac{16\epsilon^4}{1-2\epsilon}.
\end{displaymath}

\paragraph{Conclusion.}
The sequence $(\tilde{j}_t-\tilde{i}_t)_{t \geq t_0}$ is a random walk with an asymptotic drift
\begin{displaymath}
  D = \lim_{t \to \infty} \esp{\tilde{J}_t} - \lim_{t \to \infty} \esp{\tilde{I}_t} \geq 15 \epsilon^2 + 12 \epsilon^3 + \frac{32\epsilon^4}{1-2\epsilon} > 0
\end{displaymath}
when $\epsilon \in (0,1/2)$, and with bounded negative increments. So it is transient and the CA $1000$ with error $\epsilon$ is ergodic. \qed

\section{Perpectives}
In this article, the idea of studying the sizes of decorrelated islands was applied in the simple case of PCA with two-size neighbourhood and two-size alphabet, and the improvements on Section~\ref{sec:proof1000} only to the CA $1000$ with errors. These ideas could be extended and applied to a more general context such that PCA with larger size neighbourhood, with larger size alphabet and, as in Section~\ref{sec:proof1000}, we can keep more informations (the knowledge of two, three, or more boundary states) that just one boundary state as in Section~\ref{sec:ergo}. Hence, we could obtain numerous classes of ergodic general PCA with positive rates. In particular, it should be very useful to prove ergodicity of CA with errors when the CA has some stable periodic sequences of states. Nevertheless, it is impossible to prove the ergodicity of all positive rate PCA with two-size neighbourhood and two-size alphabet with just this idea. Moreover, larger the size of the neighbourhood or the alphabet is, smaller the volume of the set of ergodic PCA obtained by the idea should be in proportion of the total volume of the set of PCA.

\section{Acknowledgement}
I am very grateful to Irène Marcovici. Her knowledge on the subject confirms me that the idea developed in the paper should be written and published. Moreover, her attentive reading of the paper permits to clarify and improve it.

\bibliographystyle{alpha}
\bibliography{aaa}
\end{document}